\documentclass{ifacconf}  

\usepackage{overpic}
\usepackage{amsfonts}
\usepackage{graphicx}
\usepackage{amsmath}
\usepackage{url}
\usepackage{natbib}
\usepackage{ifthen}
\usepackage{color}

\def\eq#1{\begin{equation}#1\end{equation}}
\newcommand{\R}{{\rm I\!R}}

\def\rep#1{(\ref{#1})}
\def\scr#1{{\cal #1}}
\newtheorem{theorem}{Theorem}
\newtheorem{lemma}{Lemma}
\newtheorem{proposition}{Proposition}
\newtheorem{corollary}{Corollary}

\newtheorem{remark}{Remark}

\newcommand{\bbb}{\mathbb}

\newcommand{\1}{\mathbf{1}}
\newcommand{\0}{\mathbf{0}}

\begin{document}
\begin{frontmatter}

\title{Controlling a Networked SIS Model via a Single Input over Undirected Graphs\thanksref{footnoteinfo}}

\thanks[footnoteinfo]{This work was supported in parts by the Research Grants Council of Hong Kong Special Administrative Region, China, under the Theme-Based Research Scheme T23-701/14-N, the National Science Foundation, grants CPS 1544953 and ECCS 1509302, and the US Army Research Office (ARO) Grant W911NF-16-1-0485.}

\author[First]{Dan Wang}
\author[Second]{Ji Liu}
\author[Third]{Philip E. Par\'{e}}
\author[Fourth]{Wei Chen}
\author[First]{Li Qiu}
\author[Fifth]{Carolyn L. Beck}
\author[Fifth]{Tamer Ba\c{s}ar}

\address[First]{Department of Electronic \& Computer Engineering, The Hong Kong University of Science and Technology, Clear Water Bay, Kowloon, Hong Kong, China (e-mail: dwangah@connect.ust.hk, eeqiu@ust.hk).}
\address[Second]{Department of Electrical and Computer Engineering, Stony Brook University, NY 11794 USA (e-mail: ji.liu@stonybrook.edu)}
\address[Third]{Division of Decision and Control Systems, KTH Royal Institute of Technology, Stockholm, Sweden (e-mail: philipar@kth.se)}
\address[Fourth]{Department of Mechanics and Engineering Science \& Beijing Innovation Center for Engineering Science and Advanced Technology, Peking University, Beijing, China (e-mail: w.chen@pku.edu.cn).}
\address[Fifth]{Coordinated Science Laboratory, University of Illinois at Urbana-Champaign, IL 61801 USA (e-mail: beck3@illinois.edu, basar1@illinois.edu)}

\begin{abstract}
This paper formulates and studies the problem of controlling a networked SIS model using a single input in which the network structure is described by a connected undirected graph. A necessary and sufficient condition on the values of curing and infection rates for the healthy state to be exponentially stable is obtained via the analysis of signed Laplacians when the control input is the curing budget of a single agent. In the case when the healthy state is stabilizable, an explicit expression for the minimum curing budget is provided.
The utility of the algorithm is demonstrated using a simulation over a network of cities in the northeastern United States.

\end{abstract}

\end{frontmatter}

\section{Introduction}

Mathematical models of virus spread have been studied for centuries \citep{bernoulli1760essai}. These models have offered interesting insight into how spread processes appear in real world systems \citep{pare2018dtjournal}. This paper focuses on SIS (susceptible-infected-susceptible) models that were first introduced by \cite{kermack1932contributions}. The idea is that every member of the population is either healthy (susceptible) or infected, and transitions occur from being infected to being healthy according to some curing rate and from being healthy to infected according to some infection rate.

The original SIS model \citep{kermack1932contributions} had assumed that the population is completely connected and has instantaneous mixing, modeling the population as two groups, a susceptible group and an infected group. Researchers have extended the model to capture more realistic structures, including nontrivial graphs in the mathematical models \citep{yorke,FallMMNP07,van2009virus,KhanaferAutomatica14}. There are two interpretations of these models and their states: 1) that each node is a subpopulation and the state is a proportion of infected individuals in the corresponding subpopulation \citep{yorke,FallMMNP07}, or 2) that each node is a single individual and the state is the probability of that individual being infected \citep{van2009virus}. Various extensions of these models have been proposed to include time--varying graph structures \citep{pare2015stability,pare2016epidemic}, and competing viruses on layered networks \citep{santos2015bivirus,watkins2016optimal,bivirus_tac}.

There have been various control techniques proposed for networked SIS models; for example, \cite{stu1,wan2007network,KhanaferAutomatica14,PreciadoTCNS14,watkins2016optimal,suppress,Drakopoulos2014}, to name a few.
\cite{Drakopoulos2014} 
proposed 
a dynamic curing policy 
for containment of a contagion process modeled by a continuous Markov process. It was shown that if the CutWidth of the underlying graph is sublinear in the number of nodes, the proposed policy achieved sublinear expected time to extinction. 
\cite{Borgs2010}
studied
the epidemic threshold 
when the antidote was distributed non-uniformly. 
\cite{wan2007network,wan2008designing}
proposed and implemented 
distributed control techniques for setting healing rates and quarantine protocols 
on a severe acute respiratory syndrome (SARS) simulation model. 
The work of \cite{KhanaferAutomatica14} has proposed an antidote control technique. 
\cite{PreciadoTCNS14} and \cite{watkins2016optimal}
developed 
an optimal vaccination control technique 
using geometric programming ideas. 
\cite{suppress}
use 
distributed optimization 
to solve several formulations of curing resource
allocation problems. The proposed algorithms 
by 
\cite{suppress}, while efficient and mostly distributed, 
all appeal to an additional, independent consensus process to calculate a piece of centralized information, and thus require a synchronous stopping time across the network.
None of the other proposed algorithms in the literature are distributed. 
In fact, it has been shown 
by
\cite{bivirus_tac} that a certain type of distributed feedback controller can never stabilize the healthy state. 

With the preceding discussions
in mind, we are interested in the following 
minimum control budget problem. Suppose that the networked SIS system is unstable at the healthy state, which implies that the network will ultimately converge to an epidemic state as long as at least one agent gets infected \citep{yorke,FallMMNP07,van2009virus}, and that we can select a fixed set of agents to boost their curing rates. This technique gives a fixed, additive control input to their curing rates that is in place for all time. 
Then, is it possible to stabilize the healthy state? If possible, what is the minimum total control budget needed to stabilize the healthy state?
These problems turn out to be challenging and, to the best of our knowledge, there is no existing algorithm or analytic solution.

In this paper, we study a special case of the problems identified above, in which only one control input is allowed. 
This special case is a good starting point for several reasons. Consider the interpretation of the model where each node is a subpopulation; if a government has enough resources to treat only one subpopulation, this approach can show if the disease can be eradicated given that only one subpopulation can be treated, and if so, which is the best subpopulation to treat, and how much effort will be required to eradicate the disease.
An analogous interpretation can also be made for a computer network to be protected by equipping virus patches on only one computer.
In the case when the underlying graph is undirected, we obtain a necessary and sufficient condition for the healthy state to be exponentially stable via the mathematical tool of signed Laplacians. When the condition is satisfied, i.e., when the healthy state is stabilizable, we provide an explicit expression for the minimum control budget.


{\em Notation:}
For any positive integer $n$, we use $[n]$ to denote the set $\{1,2,\ldots,n\}$.
We use $\0$ and $\1$ to denote the vectors whose entries all equal $0$ and $1$, respectively,
and $I$ to denote the identity matrix,
while the dimensions of the vectors and matrices are to be understood from the context.
For any two real vectors $a,b\in\R^n$, we write $a\geq b$ if
$a_{i}\geq b_{i}$ for all $i\in[n]$,
$a>b$ if $a\geq b$ and $a\neq b$, and $a \gg b$ if $a_{i}> b_{i}$ for all $i\in[n]$.
We use $M^\dagger$ to denote the Moore-Penrose inverse of a matrix $M$.
The eigenvalues of a symmetric matrix $A\in\mathbb{R}^{n\times n}$ are denoted by $\lambda_1(A)\leq\dots\leq\lambda_n(A)$. 
We write $A\succ 0$ and $A\succeq 0$ when the matrix $A$ is (symmetric) positive definite and positive semidefinite, respectively. 



\section{Problem Formulation} \label{problem}

The continuous-time single-virus SIS model over an $n$-agent network (after mean-field approximation) is as follows~\citep{yorke}:
\begin{eqnarray}
\dot x_i(t) &=& -\delta_i x_i(t)+(1-x_i(t))\sum_{j=1}^n \beta_{ij}x_j(t), \nonumber\\
&&x_i(0)\in[0,1],\;\;\;\;\; i\in[n],\label{update1}
\end{eqnarray}
where $x_i(t)$ denotes the probability of agent $i$ to be infected, $\delta_i$ is the nonnegative curing rate of agent $i$, and $\beta_{ij}$ is the nonnegative infection rate between agents $i$ and $j$, whose value equals zero if they are not neighbors. The neighbor relationships among the $n$ agents are described by an undirected $n$-vertex graph $\bbb G$ whose vertices represent the agents and edges depict the neighbor relationships. We assume that $\beta_{ij} = \beta_{ji}$ for all $i,j\in[n]$.

The $n$ equations in~\rep{update1} can be combined into one equation in a state-space form as
\eq{\dot x(t) = \left(-\Delta+B-X(t)B \right)x(t), \label{sys1}}
where $x(t)$ is the state vector in $\R^n$ whose $i$th entry is $x_i(t)$,
$\Delta$ is the $n\times n$ diagonal matrix whose
$i$th diagonal entry is $\delta_i$, $B$ is the $n\times n$ matrix whose
$ij$th entry is $\beta_{ij}$, and $X(t)$ is the $n\times n$ diagonal matrix
whose $i$th diagonal entry is $x_i(t)$.
Since we have assumed that $\beta_{ij} = \beta_{ji}$ for all $i,j\in[n]$, $B$ is symmetric. 
It is also assumed that $B$ is an irreducible matrix, which implies that the underlying graph $\bbb G$ is 
connected.

It is easy to see that $\0$ is an equilibrium of system~\rep{sys1}. We call it the healthy state. 
It has been shown in \cite{yorke,van2009virus,FallMMNP07} that the healthy state is the unique equilibrium of the system if $\Delta-B \succeq 0$, and if $\Delta-B \succeq 0$ does not hold, the healthy state is unstable, and there will be a unique epidemic equilibrium $x^*\gg \0$, which is almost globally stable. More can be said.
The following result is a variant of Proposition III.1 in~\cite{preciado2017optimal}.


\begin{proposition}
The healthy state $x=\0$ of system~\rep{sys1} is exponentially stable if, and only if, $\Delta-B \succ 0$. 
\label{stability}\end{proposition}

With this in mind, we formulate the following control problem which aims to stabilize the healthy state using a single control input.


Suppose that the necessary and sufficient condition in Proposition \ref{stability} is not met, and we have a control budget to increase the curing rate of a single agent $i$, i.e, 
\begin{eqnarray}
\dot x_i(t) &=& -(\delta_i+u_i) x_i(t)+(1-x_i(t))\sum_{j=1}^n \beta_{ij}x_j(t), \nonumber\\
\dot x_k(t) &=& -\delta_k x_k(t)+(1-x_k(t))\sum_{j=1}^n \beta_{kj}x_j(t), \nonumber\\
&&k\neq i,\;\;\;\;\; k\in[n].\label{problem1}
\end{eqnarray}
We seek to answer the questions: what is the necessary and sufficient condition to stabilize the healthy state? When the healthy state is stabilizable with such a single control input, what is the minimum control budget? 
Formally stated we have:
\begin{equation}
\begin{aligned}
& \underset{i\in[n]}{\text{minimize}}
& & u_i \\
& \text{subject to}& & \Delta + U_i-B \succ 0, 
\end{aligned}
\label{p}\end{equation}
where $U_i$ is the $n\times n$ diagonal matrix whose $i$th diagonal entry equals $u_i$ and the other diagonal entries all equal zero.
In the following, we will present machinery that enables solving this problem in an optimal and efficient manner.

Mathematically, the problems in (3) and (4) amount to making a symmetric matrix positive definite by changing one of its diagonal entries with minimal effort. We are aware that this problem, in its more general form of changing multiple diagonal entries, has been explored in the literature based on linear matrix inequalities (LMI) \citep{BhatiaBook}, which provides numerical solutions but not analytical ones. 

The recent works \citep{stu1,stu2} studied the problem of minimizing the dominant eigenvalue of $\Delta + U_i-B$ when a subset of nodes can be controlled, subject to a given total control budget, with $u_i\geq 0$ in the first paper and $u_i\in \mathbb{R}$ in the second one. Algorithms have been proposed for optimal subset design. This problem is different from the one we consider in this paper. Here, we wish to minimize the total control budget such that $\Delta + U_i-B\succ 0$. 
Moreover, we want to obtain an analytic solution for the minimum control budget when the healthy state is stabilizable. We achieve these goals for the special case of changing a single diagonal entry by exploiting the structure of signed graphs.



\section{Main Results} \label{result}

Let $D_i$ be the $n\times n$ diagonal matrix whose diagonal entries are 
\begin{eqnarray*}
d_i &=& 0, \\
d_k &=& \delta_k-\sum_{j=1}^n\beta_{kj}, \;\;\; k\neq i, \;\;\; k\in[n]. 
\end{eqnarray*}
Let $V_i$ be an $n\times n$ diagonal matrix whose $i$th diagonal entry equals $(\delta_i-\sum_{j=1}^n\beta_{ij})$ and all the other diagonal entries equal zero. Define
\begin{equation}
L_i=\left[\begin{array}{cc}
               \Delta-V_i-B & -D_i\mathbf{1}\\
               -\mathbf{1}^\top D_i&\mathbf{1}^\top D_i\mathbf{1}\end{array}\right].
\label{phil}\end{equation}
Note that $\mathbf{1}^\top D_i\mathbf{1} = \sum_{j=1}^n d_j$. 
By construction, the
row sums of $L_i$ all equal zero.

The following theorem provides an explicit lower bound on the control input $u_i$ that stabilizes the healthy state. 

\begin{theorem}\label{thm}
Suppose that $u_i$ is the single control input of agent $i$.
The healthy state of system~\rep{problem1} is exponentially stable if, and only if, $L_i\succeq 0$ with a simple eigenvalue at zero and 
\begin{align}
u_i > \sum_{j=1}^n\beta_{ij} - \delta_i - 
\left((e_i-e_{n+1})^\top L_{i}^{\dagger}(e_i-e_{n+1})\right)^{-1}, \label{bound}
\end{align} 
where $e_i$ is the $(n+1)$-dimensional unit vector with the $i$th entry being one and others zero.
\end{theorem}

It is worth emphasizing that $L_i$ is independent of $u_i$, as is the right side of~\rep{bound}. The matrix $L_i\succeq 0$ and $0$ is a simple eigenvalue with eigenspace $\mathrm{span}\{\1\}$. Thus $L_i^\dagger\succeq 0$ and $0$ is a simple eigenvalue with the eigenspace $\mathrm{span}\{\1\}$, which will be shown in the next section. Since $(e_i-e_{n+1})\notin \mathrm{span}\{\1\}$, the last term in (\ref{bound}) is positive. If the value of the right side of~\rep{bound} is negative, it implies that the healthy state is exponentially stable without any control input, i.e., $\Delta-B \succ 0$.

Given any $i$, by Theorem~\ref{thm}, since $L_i$ is independent of $u_i$, 
we can first check whether the curing resources on the other agents are sufficient to stabilize the healthy state or not, i.e., whether the condition that $L_i\succeq 0$ with a simple eigenvalue at zero holds or not. Subsequently, if that condition is met, then we can design a control input, appealing to \eqref{bound}, to make the healthy state exponentially stable.


\begin{corollary}\label{main}
Suppose that $u$ is the single control input on any one agent.
The healthy state of system~\rep{problem1} is stabilizable if, and only if, there exists at least one $i\in[n]$ such that $L_i\succeq 0$ with a simple eigenvalue at zero. 
Let $\scr J$ be the set of such indices in $[n]$.  
If the condition holds, i.e., $\scr J \neq \emptyset$, the minimum control input to stabilize the healthy state is given by 
$$u=\inf_{i\in\scr J} \; u_i,$$
where 
$u_i$ satisfies~\rep{bound}. 
\end{corollary}



The corollary is a direct consequence of Theorem~\ref{thm} and provides the machinery to explicitly solve the problem posed in~\rep{p}. 

\section{Analysis and Proofs} \label{analysis}

To prove the main results stated in the previous section, we need the following preliminaries on results for Laplacian matrices and their variants.

\subsection{Laplacian Matrices}

Consider a weighted undirected graph $\mathbb{G}=(\mathcal{V},\mathcal{E})$ with no self-loops, which consists of a set of nodes $\mathcal{V}=[n]$ and a set of edges $\mathcal{E}$. 
We associate with each edge $(i,j)\in\mathcal{E}$ a real-valued weight $w_{ij}$. The associated Laplacian matrix $L=[l_{ij}]_{n\times n}$ is defined by
\begin{align}\label{laplacian}
\begin{split}
l_{ij}=\begin{cases}-w_{ij},&i\neq j,\\ \sum_{k\neq i} w_{ik},& i=j.\end{cases}
\end{split}
\end{align}
It is clear that $L$ is a symmetric matrix whose row sums all equal zero, and thus zero is an eigenvalue of $L$.
It is worth emphasizing that we allow negative weights $w_{ij}$. 
In the case where all weights $w_{ij}$ are nonnegative, it is well known that $L$ is positive semidefinite. Moreover, zero is a simple eigenvalue with eigenspace 
$\mathrm{span}\{\mathbf{1}\}$ if and only if $\bbb G$ is connected \citep{fiedler}.  
In the case when there exists at least one negative weight $w_{ij}$, $L$ is sometimes called a signed Laplacian, whose property of positive semidefiniteness has been recently investigated in \cite{zelazo2014}, \cite{yongxin}, and \cite{cdc16}. 
With this work in mind, any symmetric matrix whose row sums all equal zero can be regarded as a Laplacian matrix whose corresponding weighted undirected graph may contain negative weights.

There is a useful factorization of a Laplacian matrix $L$. 
Suppose that the corresponding graph $\bbb G$ has $m$ edges. 
Let
$W\!=\!\mathrm{diag}\{w_1,w_2,\dots,w_m\}$ denote the $m\times m$ diagonal matrix where the diagonal entries are the edge weights, i.e., for each edge $\varepsilon_k\in\mathcal{E}$,
$w_{k}=w_{ij}$ for $(i,j)= \varepsilon_k$.
The order of the $m$ edges can be arbitrary. 
In addition, we assign an arbitrary orientation to each edge, i.e., for each edge $\varepsilon_k\in\mathcal{E}$, we set an arbitrary endpoint as the head and the
other one as the tail. The oriented incidence matrix $B=[b_{ij}]\in\mathbb{R}^{n\times m}$ is defined as
\begin{align*}
b_{ik}=\begin{cases}1, &\text{if $i$ is the head of $\varepsilon_k$},\\-1, &\text{if $i$ is the tail of $\varepsilon_k$},\\0,&\text{otherwise}.\end{cases}
\end{align*}
Then, the Laplacian matrix $L$ can be factorized as
\begin{equation}
L=BWB^\top.
\label{incidence}\end{equation}
It is worth noting that while the incidence matrix $B$ depends on the choice of orientations, the Laplacian matrix $L$ does not.
In the case when both positive and negative weights exist, we define $\bbb{G}_+$ and $\bbb{G}_-$ as the spanning subgraphs of $\bbb G$ which consist of all positive- and negative-weighted edges, and we call them positive and negative subgraphs, respectively. In this case, $L$ can also be written in the following form:
\begin{equation}
L=B_+W_+B_+^\top + B_-W_-B_-^\top,
\label{decomp}\end{equation}
where $B_+$ and $B_-$ are incidence matrices corresponding
to the positive and negative subgraphs, respectively, and $W_+$ and $W_-$ are the diagonal matrices whose diagonal entries are the weights of the edges in the positive and negative subgraphs, respectively.
From~\rep{incidence}, it is easy to see that $B_+W_+B_+^\top$ and $B_-W_-B_-^\top$ are the Laplacian matrices of $\bbb{G}_+$ and $\bbb{G}_-$, respectively.
In fact, it can be verified that if $\bbb G$ can be partitioned into $p$ spanning subgraphs, $\bbb G_1,\bbb G_2,\ldots,\bbb G_p$, where all $p$ subgraphs have the same set of vertices, if their edge sets are mutually disjoint, and if the union of their edge sets equals the edge set of $\bbb G$, then the Laplacian $L$ can be written as 
\begin{equation}
L=\sum_{i=1}^p B_iW_iB_i^\top,
\label{partition}\end{equation}
where $B_i$ is the incidence matrix of $\bbb G_i$ and $W_i$ is the diagonal matrix whose diagonal entries are the weights of the edges in $\bbb G_i$.
It can be seen that each $B_iW_iB_i^\top$ equals the Laplacian matrix of $\bbb{G}_i$.

We will later make use of the following lemmas.

\begin{lemma}
A Laplacian matrix $L$ satisfies $L\succeq 0$ with a simple eigenvalue at zero and associated eigenspace $\mathrm{span}\{\mathbf{1}\}$ if, and only if, $L^{\dagger}\succeq 0$ with a simple eigenvalue at zero and associated eigenspace $\mathrm{span}\{\mathbf{1}\}$.
\label{xxx}\end{lemma}

{\em Proof:} 
Consider the Schur decomposition of $L$ as
\[
L=\begin{bmatrix}U & \frac{1}{\sqrt{n}}\1 \end{bmatrix} \begin{bmatrix}
S & 0 \\0 & 0\end{bmatrix} \begin{bmatrix}U & \frac{1}{\sqrt{n}}\1 \end{bmatrix}^\top,
\]
where $\begin{bmatrix}U & 1/{\sqrt{n}}\1 \end{bmatrix}$ is a unitary matrix and $S$ is a diagonal matrix containing all the nonzero eigenvalues of $L$. By \citep[Theorem 1.2.1]{CampbellMeyerBook}, we have
\begin{align*}
L^\dagger&=\begin{bmatrix}U & \frac{1}{\sqrt{n}}\1 \end{bmatrix} \begin{bmatrix}
S & 0 \\0 & 0\end{bmatrix}^\dagger \begin{bmatrix}U & \frac{1}{\sqrt{n}}\1 \end{bmatrix}^\top \\
&=\begin{bmatrix}U & \frac{1}{\sqrt{n}}\1 \end{bmatrix} \begin{bmatrix}
S^{-1} & 0 \\0 & 0\end{bmatrix} \begin{bmatrix}U & \frac{1}{\sqrt{n}}\1 \end{bmatrix}^\top,
\end{align*}
which is the Schur decomposition of $L^\dagger$. Therefore, the statement holds.

\begin{lemma}
Let $L$ be an $n\times n$ Laplacian matrix and $Q$ be a matrix in $\mathbb{R}^{n\times (n-1)}$ such that
$Q^\top Q = I$ and $Q^\top \1=0$.
Then, $L\succeq 0$ with simple zero eigenvalue if, and only if, $Q^\top LQ\succ 0$.
\label{proj}\end{lemma}

{\em Proof:}
Since $Q^\top Q = I$ and $Q^\top \1=0$, it follows that 
$$\begin{bmatrix}Q & \frac{1}{\sqrt{n}}\1 \end{bmatrix}$$
is a unitary matrix. Note that 
\begin{align*}
\begin{bmatrix}Q & \frac{1}{\sqrt{n}}\1 \end{bmatrix}^\top L \begin{bmatrix}Q & \frac{1}{\sqrt{n}}\1 \end{bmatrix}
=\begin{bmatrix}
Q^\top LQ & 0 \\0 & 0 \end{bmatrix}.
\end{align*}
Then, the spectrum of $L$ consists of zero and the eigenvalues of $Q^\top LQ$, from which the lemma follows.

\begin{lemma}
Let $L\succeq 0$ be an $n\times n$ Laplacian matrix with a simple zero eigenvalue and $Q$ be a matrix in $\mathbb{R}^{n\times (n-1)}$ such that
$Q^\top Q = I$ and $Q^\top \1=0$.
Then, $L^{\dagger}=Q(Q^\top L Q)^{-1}Q^\top$.
\label{yyy}\end{lemma}

{\em Proof:} 
It can be seen that
$
\begin{bmatrix}Q & 1/{\sqrt{n}}\1 \end{bmatrix}
$
is a unitary matrix. By Theorem 1.2.1 in \cite{CampbellMeyerBook} and Lemma~\ref{proj}, we have
\begin{align*}
L^{\dagger}&=\begin{bmatrix}Q & \frac{1}{\sqrt{n}}\1 \end{bmatrix} \begin{bmatrix}Q & \frac{1}{\sqrt{n}}\1 \end{bmatrix}^\top L^{\dagger} \begin{bmatrix}Q & \frac{1}{\sqrt{n}}\1 \end{bmatrix} \begin{bmatrix}Q & \frac{1}{\sqrt{n}}\1 \end{bmatrix}^\top\\
&=\begin{bmatrix}Q & \frac{1}{\sqrt{n}}\1\end{bmatrix} \left(\begin{bmatrix}Q & \frac{1}{\sqrt{n}}\1 \end{bmatrix}^\top L \begin{bmatrix}Q & \frac{1}{\sqrt{n}}\1 \end{bmatrix}\right)^{\dagger} \begin{bmatrix}Q & \frac{1}{\sqrt{n}}\1\end{bmatrix}^\top\\
&=\begin{bmatrix}Q & \frac{1}{\sqrt{n}}\1\end{bmatrix} \begin{bmatrix}
Q^\top LQ & 0 \\ 0 & 0\end{bmatrix}^{\dagger} \begin{bmatrix}Q & \frac{1}{\sqrt{n}}\1\end{bmatrix}^\top \\
&=\begin{bmatrix}Q & \frac{1}{\sqrt{n}}\1\end{bmatrix} \begin{bmatrix}
(Q^\top LQ)^{-1} & 0 \\ 0 & 0\end{bmatrix} \begin{bmatrix}Q & \frac{1}{\sqrt{n}}\1\end{bmatrix}^\top \\
&=Q(Q^\top L Q)^{-1}Q^\top,
\end{align*}
which completes the proof.

\subsection{Loopy Laplacian} \label{sec:loopy}
Consider a weighted undirected graph $\mathbb{G}=(\mathcal{V},\mathcal{E})$ which consists of a set of nodes $\mathcal{V}=[n]$ and a set of edges including self-loops. 
Associate with each edge $(i,j)\in\mathcal{E}$ a real-valued weight $w_{ij}$. The weight of the self-loop of node $i$ is denoted by $w_{ii}\in\mathbb{R}$. 
The associated loopy Laplacian $H$ \citep{Dorfler13} is defined as
\begin{equation*}
H=L+D,
\end{equation*}
where $L$ is the Laplacian matrix as in~\rep{laplacian}, and
$D=\textrm{diag}\{w_{11},\dots,w_{nn}\}$ is a diagonal matrix whose diagonal entries are the self-loop weights. When $D\neq 0$, i.e., at least one of $w_{ii}$ is nonzero, $H$ is called a strictly loopy Laplacian \citep{Dorfler13}. Clearly, $H$ and $L$ are symmetric matrices.

Given a strictly loopy Laplacian, we introduce a grounded node with index $n+1$, and connect node $i\in[n]$ to the grounded node when $w_{ii}\neq 0$; the self-loops are then eliminated. This process gives us a loopless graph $\hat{\mathbb{G}}$ (i.e., there are no self-loops in the graph) with either positive or negative weights, and its signed Laplacian $\hat{L}\in\mathbb{R}^{(n+1)\times(n+1)}$. The augmented Laplacian $\hat{L}$ is defined as
\begin{equation}
\hat{L}=\left[\begin{array}{cc}
               H&-D\mathbf{1}\\
               -\mathbf{1}^\top D&\sum_{i=1}^n w_{ii}\end{array}\right],
\label{augment}
\end{equation}
which is also a symmetric matrix.

\begin{lemma}
\label{lemma2}
Consider a strictly loopy Laplacian matrix $H\in\mathbb{R}^{n\times n}$ and its corresponding augmented Laplacian matrix $\hat{L}$. Then, $H\succ 0$ if, and only if, $\hat{L}\succeq 0$ with a simple zero eigenvalue.
\end{lemma}

{\em Proof:}
We first prove the sufficiency. Suppose that $\hat{L}\succeq 0$ and zero is a simple eigenvalue. 
Since $\hat{L}\mathbf{1}=0$ and zero is a simple eigenvalue, the eigenspace of the zero eigenvalue is $\mathrm{span}\{\mathbf{1}\}$.
Let $x\in\mathbb{R}^n$ be any nonzero vector and set $y=[x^\top\ 0]^\top$. Since $y\notin\mathrm{span}\{\mathbf{1}\}$, from \rep{augment}, we have
$x^\top Hx=y^\top \hat{L}y>0$,
which implies that $H$ is positive definite.

We now prove the necessity. Suppose that $H\succ 0$. Denote the eigenvalues of $H$ and $\hat{L}$ by
$\lambda_1(H) \leq \dots \leq \lambda_n(H)$
and 
$\lambda_1(\hat{L}) \leq \dots \leq \lambda_{n+1}(\hat{L})$,
respectively. Since $H$ is a principal submatrix of $\hat{L}$, by the Cauchy Interlace Theorem (see Theorem~4.3.8 in \cite{Horn}), we have
$\lambda_1(\hat{L})\leq \lambda_1(H) \leq \lambda_2(\hat{L})\leq \dots \leq \lambda_n(H) \leq \lambda_{n+1}(\hat{L})$.
Thus, $\lambda_i(\hat{L})>0$ for all $i\in\{2,\dots,n+1\}$. Since the row sums of $\hat{L}$ are all equal to zero, $\hat{L}$ has a zero eigenvalue, which implies that $\lambda_1(\hat{L})=0$. This completes the proof.


\subsection{Proof of Theorem~\ref{thm}}

To prove the theorem, we will need the following lemma. 

\begin{lemma}
{\rm (See, e.g., Theorem 1.10 in \cite{dullerud2013course})}
Suppose that $M$ is a symmetric matrix which is partitioned into four submatrices as
$$
M = \left[\begin{array}{cc}
               A & B\\
               B^\top & D\end{array}\right],
$$
where $A$ and $D$ are square matrices. Then,
\begin{itemize}
\item [1)] $M\succ 0$ if, and only if, $A\succ 0$ and $D-B^\top A^{-1}B \succ 0$;
\item [2)] $M\succ 0$ if, and only if, $D\succ 0$ and $A - BD^{-1}B^\top \succ 0$.
\end{itemize}
\label{schur}\end{lemma}

To proceed, we rewrite the matrix $\Delta+U_i-B$ as follows. 
Let $\hat{D}_i$ be the diagonal matrix whose diagonal entries are 
\begin{eqnarray}
\hat{d}_i &=& \delta_i+u_i-\sum_{j=1}^n\beta_{ij}, \\
\hat{d}_k &=& \delta_k-\sum_{j=1}^n\beta_{kj}, \;\;\; k\neq i, \;\;\; k\in[n].
\end{eqnarray}
Then, $\Delta+U_i-B-\hat{D}_i$ is a matrix whose row sums all equal zero, and thus can be regarded as a Laplacian matrix. 
Thus, $\Delta+U_i-B=(\Delta+U_i-B-\hat{D}_i)+\hat{D}_i$ can be treated as a strictly loopy Laplacian of a weighted graph in which the weight on each edge $(i,j)$ is $\beta_{ij}$ and the weight of the self-loop of each node $j$ is $\hat{d}_j$. 
To proceed, let 
\begin{equation}
\hat{L}_i=\left[\begin{array}{cc}
               \Delta+U_i-B & -\hat{D}_i\mathbf{1}\\
               -\mathbf{1}^\top \hat{D}_i&\mathbf{1}^\top \hat{D}_i\mathbf{1}\end{array}\right].
\label{dan}\end{equation}
Note that $\mathbf{1}^\top \hat{D}_i\mathbf{1} = \sum_{j=1}^n \hat{d}_j$. It is easy to check that all row sums of $\hat{L}_i$ are equal to zero. 
Note, $\hat{L}_i$ is the augmented Laplacian for $\Delta+U_i-B$. 
Specifically, 
we have introduced a grounded node with index $n+1$ to the existing graph $\bbb G$, as explained in Section \ref{sec:loopy}, 
to create a loopless graph $\hat{\mathbb{G}}$, and thus removing the self-loops from the graph. 
The corresponding signed Laplacian $\hat{L}_i\in\mathbb{R}^{(n+1)\times(n+1)}$, given in~\rep{dan}, has 
possibly both positive and negative weights.

Note that the difference between $\hat{D}_i$ and $D_i$ only lies in the $i$th diagonal entry. It is not hard to verify that 
$L_i$ given in~\rep{phil} is the Laplacian matrix of the graph obtained by removing the edge $(i,n+1)$ from $\hat{\mathbb{G}}$.

From Proposition~\ref{stability} and Lemma~\ref{lemma2}, it is sufficient to show that $\hat{L}_i\succeq 0$ with a simple zero eigenvalue if and only if $L_i\succeq 0$ with a simple zero eigenvalue and \eqref{bound} holds.

Let $Q$ be a matrix in $\mathbb{R}^{(n+1)\times n}$ such that 
$$
\begin{bmatrix}Q& \frac{1}{\sqrt{n+1}}\mathbf{1}\end{bmatrix}
$$ 
is a unitary matrix. Consequently, $Q^\top Q = I$ and $Q^\top \1=0$. 
Then, from Lemma~\ref{proj}, $\hat{L}_i\succeq 0$ with a simple zero eigenvalue if and only if $Q^\top\hat{L}_iQ\succ 0$.
Since $\hat{L}_i$ is a Laplacian matrix whose corresponding graph $\hat{\bbb G}$'s edges are all positive-weighted except for the edge $(i,n+1)$ whose sign depends on the value of $u_i$. 
From~\rep{partition}, 
$$
\hat{L}_i = L_i + (e_i-e_{n+1})\hat{d}_i(e_i-e_{n+1})^\top,
$$
in which $(e_i-e_{n+1})\hat{d}_i(e_i-e_{n+1})^\top$ equals the Laplacian of the spanning subgraph of $\hat{\bbb G}$ whose edge set has only one edge $(i,n+1)$.
Thus, $\hat{L}_i\succeq 0$ with a simple zero eigenvalue if and only if
\begin{align}
Q^\top\left(L_i + (e_i-e_{n+1})\hat{d}_i(e_i-e_{n+1})^\top\right)Q\succ 0.\label{ine1}
\end{align}
If $\hat{d}_i>0$, then $\hat{\bbb G}$ is a connected graph whose edge weights are all positive, which implies that $\hat{L}_i \succeq 0$ with a simple eigenvalue at zero. 
If $\hat{d}_i=0$, then $\hat{L}_i =L_i$.
From Lemma~\ref{xxx}, it can be seen that $L_i\succeq 0$ is equivalent to \eqref{bound}, which implies that the statement  
$\hat{L}_i\succeq 0$ with a simple zero eigenvalue is equivalent to the statement $L_i\succeq 0$ with a simple zero eigenvalue and \eqref{bound} holds.
We therefore consider the case when $\hat{d}_i<0$.

\begin{figure}
\centering
\includegraphics[width = .65\columnwidth]{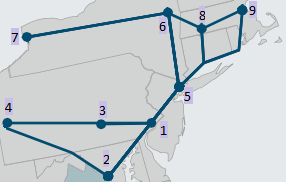}
\caption{
Plot of the northeast train routes: we ignore the routes that lead out of this subset of cities.}
\label{fig:amtrak}
\end{figure}

From item 1) in Lemma~\ref{schur}, \eqref{ine1} is equivalent to
\begin{align}
\begin{bmatrix}\left(-\hat{d}_i\right)^{-1} & (e_i-e_{n+1})^\top Q\\Q^\top (e_i-e_{n+1}) & Q^\top L_i Q\end{bmatrix}\succ 0.\label{ine2}
\end{align}
From item 2) in Lemma~\ref{schur},
\eqref{ine2} is equivalent to the following two conditions:
\begin{align}
Q^\top L_i Q \succ 0,\label{ine3}
\end{align}
\begin{align}
-\hat{d}_i^{-1}-(e_i-e_{n+1})^\top Q(Q^\top L_i Q)^{-1}Q^\top (e_i-e_{n+1}) > 0.\label{ine4}
\end{align}
The first condition \eqref{ine3} is equivalent to that $L_i\succeq 0$ with a simple zero eigenvalue because of Lemma~\ref{proj}.
The second condition \eqref{ine4} is equivalent to \eqref{bound} due to Lemma~\ref{yyy}. 
This then completes the proof.

\begin{remark}
{
A key element in the above proof is showing that the condition that $\hat{L}_i$ is positive semidefinite with a simple zero eigenvalue is equivalent to the condition that $L_i$ is positive semidefinite with a simple zero eigenvalue and \eqref{bound} holds. Alternative proofs of such equivalence are available by adapting the techniques in the proofs of Theorem III.3 in \cite{zelazo2014} and Theorem 3 in \cite{song2018}. Still, the new proof provided as above is much more concise compared to the existing techniques in the literature.
}
\end{remark}



\section{Simulation}\label{simulation}

For the simulations we consider a virus spreading  among the northeast of the United States. The graph structure is comprised of nine cities, or subpopulations, in the northeast (see Figure \ref{fig:amtrak}). 
We make the simplifying assumption that there is a virus which spreads among subpopulations, grouped by cities, according to the passenger train network (ignoring spread among air and car travel). 
The connectivity is determined by the Amtrak Rail Service between these cities. See Figure \ref{fig:amtrak} for a plot of the train routes. The corresponding infection rate 
matrix is 
\begin{equation}\footnotesize
B_{NE} = 
\begin{bmatrix}
  1 & 1 & 1 & 0 & 1 & 0 & 0 & 0 & 0 \\
  1 & 1 & 0 & 1 & 0 & 0 & 0 & 0 & 0 \\
  1 & 0 & 1 & 1 & 0 & 0 & 0 & 0 & 0 \\
  0 & 1 & 1 & 1 & 0 & 0 & 0 & 0 & 0 \\
  1 & 0 & 0 & 0 & 1 & 1 & 0 & 1 & 1 \\
  0 & 0 & 0 & 0 & 1 & 1 & 1 & 1 & 0 \\
  0 & 0 & 0 & 0 & 0 & 1 & 1 & 0 & 0 \\
  0 & 0 & 0 & 0 & 1 & 1 & 0 & 1 & 1 \\
  0 & 0 & 0 & 0 & 1 & 0 & 0 & 1 & 1
\end{bmatrix}.
\label{eq:Btrain}
\end{equation}
We set the healing rate to $\delta_i = 3.5$ for all $i$. 
The state of the system is the proportion of infected individuals in the corresponding city. 
Red indicates completely infected ($x_i=1$) and blue indicates completely healthy ($x_i=0$).

The necessary and sufficient condition from Proposition~\ref{stability} is not satisfied; the minimum eigenvalue of $\Delta - B$ is $-0.3134$. Therefore, with no control effort the system converges to the epidemic equilibrium  
$$\footnotesize x^* = \begin{bmatrix}
.063 &\!\! .034 &\!\! .034 &\! .026 &\!\! .104 &\!\! .080&\!\!  .031 &\!\! .090 &\!\! .070\end{bmatrix}^\top,$$
meaning between 3-10\% of the populations of the cities are infected with the virus at steady state.

\begin{figure}
\centering
\includegraphics[width = .4\columnwidth]{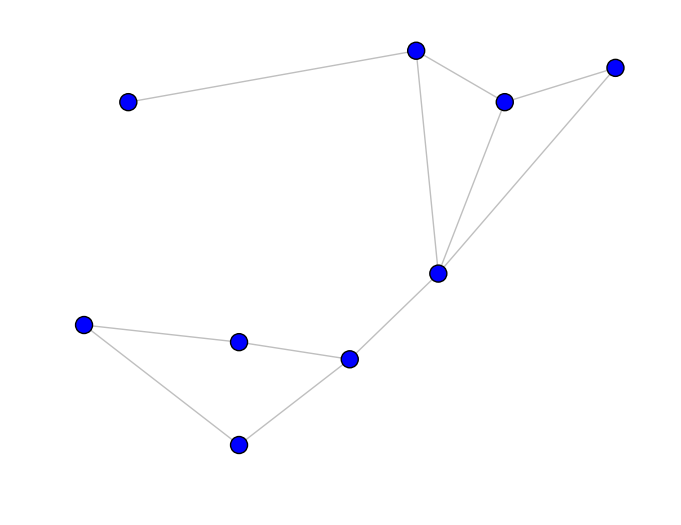}
\begin{overpic}[width = .58\columnwidth]{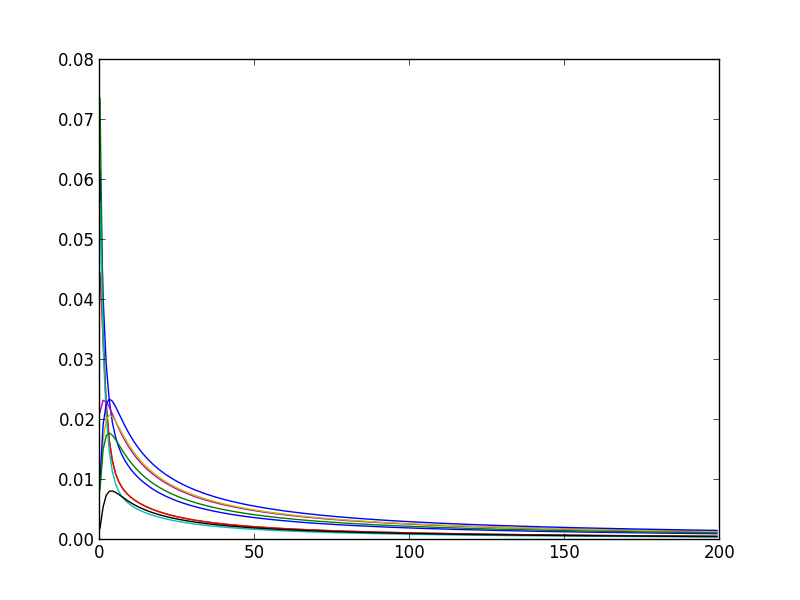}
    \put(49,-1){{\parbox{\linewidth}{%
     $t$}}}
     \put(-2,36){{\parbox{\linewidth}{%
     \rotatebox{90}{$x_i(t)$}}}}
\end{overpic}
\caption{
The healthy equilibrium of the northeast system (left, blue indicates healthy) and trajectories (right) with infection rates and connectivity given in \eqref{eq:Btrain}, the healing rate is $\delta_i = 3.5$ for all $i$, and the control input $u_5 = 1.30$. 
}
\label{fig:0}
\end{figure}

Now the question becomes the following: if we are limited by our infrastructure and can treat only one city, which one is the best one and how much treatment is required? Implementing the algorithm from Corollary \ref{main} to solve Problem \eqref{p}, answers these questions. It turns out that the only node that meets the base requirement of $L_i\succeq 0$ with a simple eigenvalue at zero is $i=5$, that is, New York City (node 5). The minimum amount of control input that eradicates the virus is 
$u_5 = 1.30$.
Including this control input in the system as in \eqref{problem1} results in the system converging to the healthy state, depicted in Figure~\ref{fig:0}.

From this example, we
discovered, as would be expected since there is only one control input, that the virus cannot be too strong for a system to be stabilizable using this approach. If $\delta_i$ is much smaller than $3.5$ in the system considered here, then there is not an $i$ such that $L_i\succeq 0$. 

%

\section{Conclusion}\label{con}
We have proposed a control technique that employs one control input on a single node to mitigate the spread of a virus in a network, that can drive the system to the healthy state. A necessary and sufficient condition on the values of curing and infection rates for the healthy state to be exponentially stable has been obtained. In the case when the healthy state is stabilizable, an explicit expression for the minimum curing budget has been provided. The derivation of the controller used a loopy Laplacian approach, employing a new, more concise proof technique than what can be found in the literature. We have shown the utility of the proposed controller via simulation.

For future work we would like to explore, using a similar control technique, situations with more than one control input. We would also like to extend the technique to be distributed across the network. Finally, we would like to implement these techniques in some real spread applications, to eradicate virus from a real system.

\bibliography{bib,hermitian}

\begin{thebibliography}{31}
\providecommand{\natexlab}[1]{#1}
\providecommand{\url}[1]{\texttt{#1}}
\providecommand{\urlprefix}{URL }
\expandafter\ifx\csname urlstyle\endcsname\relax
  \providecommand{\doi}[1]{doi:\discretionary{}{}{}#1}\else
  \providecommand{\doi}{doi:\discretionary{}{}{}\begingroup
  \urlstyle{rm}\Url}\fi

\bibitem[{Bernoulli(1760)}]{bernoulli1760essai}
Bernoulli, D. (1760).
\newblock Essai d’une nouvelle analyse de la mortalit{\'e} caus{\'e}e par la
  petite v{\'e}role et des avantages de l’inoculation pour la pr{\'e}venir.
\newblock \emph{Histoire de l’Acad. Roy. Sci. avec M{\'e}m. des Math. et
  Phys. and M{\'e}m}.

\bibitem[{Bhatia(2009)}]{BhatiaBook}
Bhatia, R. (2009).
\newblock \emph{Positive Definite Matrices}.
\newblock Princeton University Press.

\bibitem[{Borgs et~al.(2010)Borgs, Chayes, Ganesh, and Saberi}]{Borgs2010}
Borgs, C., Chayes, J., Ganesh, A., and Saberi, A. (2010).
\newblock How to distribute antidote to control epidemics.
\newblock \emph{Random Structures \& Algorithms}, 37(2), 204--222.

\bibitem[{Campbell and Meyer(2009)}]{CampbellMeyerBook}
Campbell, S.L. and Meyer, C.D. (2009).
\newblock \emph{Generalized Inverses of Linear Transformations}.
\newblock SIAM.

\bibitem[{Chen et~al.(2016{\natexlab{a}})Chen, Liu, Chen, Khong, Wang,
  Ba{\c{s}}ar, Qiu, and Johansson}]{cdc16}
Chen, W., Liu, J., Chen, Y., Khong, S.Z., Wang, D., Ba{\c{s}}ar, T., Qiu, L.,
  and Johansson, K.H. (2016{\natexlab{a}}).
\newblock Characterizing the positive semidefiniteness of signed {L}aplacians
  via effective resistances.
\newblock In \emph{Proc. 55th IEEE Conf. Decision Control}, 985--990.

\bibitem[{Chen et~al.(2016{\natexlab{b}})Chen, Khong, and Georgiou}]{yongxin}
Chen, Y., Khong, S.Z., and Georgiou, T.T. (2016{\natexlab{b}}).
\newblock On the definiteness of graph {L}aplacians with negative weights:
  geometrical and passivity-based approaches.
\newblock In \emph{Proc. American Control Conference}, 2488--2493.

\bibitem[{Dorfler and Bullo(2013)}]{Dorfler13}
Dorfler, F. and Bullo, F. (2013).
\newblock Kron reduction of graphs with applications to electrical networks.
\newblock \emph{IEEE Transactions on Circuits and Systems {I}: Regular Papers},
  60(1), 150--163.

\bibitem[{Drakopoulos et~al.(2014)Drakopoulos, Ozdaglar, and
  Tsitsiklis}]{Drakopoulos2014}
Drakopoulos, K., Ozdaglar, A., and Tsitsiklis, J.N. (2014).
\newblock An efficient curing policy for epidemics on graphs.
\newblock \emph{IEEE Transactions on Network Science and Engineering}, 1(2),
  67--75.

\bibitem[{Dullerud and Paganini(2013)}]{dullerud2013course}
Dullerud, G.E. and Paganini, F. (2013).
\newblock \emph{A Course in Robust Control Theory: A Convex Approach}.
\newblock Springer.

\bibitem[{Fall et~al.(2007)Fall, Iggidr, Sallet, and Tewa}]{FallMMNP07}
Fall, A., Iggidr, A., Sallet, G., and Tewa, J.J. (2007).
\newblock Epidemiological models and {L}yapunov functions.
\newblock \emph{Mathematical Modelling of Natural Phenomena}, 2(1), 62--83.

\bibitem[{Fiedler(1973)}]{fiedler}
Fiedler, M. (1973).
\newblock Algebraic connectivity of graphs.
\newblock \emph{Czech. Math. J.}, 23(98), 298--305.

\bibitem[{Horn and Johnson(1990)}]{Horn}
Horn, R.A. and Johnson, C.R. (1990).
\newblock \emph{Matrix Analysis}.
\newblock Cambridge University Press.

\bibitem[{Kermack and McKendrick(1932)}]{kermack1932contributions}
Kermack, W.O. and McKendrick, A.G. (1932).
\newblock Contributions to the mathematical theory of epidemics. {II}. {T}he
  problem of endemicity.
\newblock \emph{Proceedings of the Royal Society A}, 138(834), 55--83.

\bibitem[{Khanafer et~al.(2016)Khanafer, Ba\c{s}ar, and
  Gharesifard}]{KhanaferAutomatica14}
Khanafer, A., Ba\c{s}ar, T., and Gharesifard, B. (2016).
\newblock Stability of epidemic models over directed graphs: a positive systems
  approach.
\newblock \emph{Automatica}, 74, 126--134.

\bibitem[{Lajmanovich and Yorke(1976)}]{yorke}
Lajmanovich, A. and Yorke, J.A. (1976).
\newblock A deterministic model for gonorrhea in a nonhomogeneous population.
\newblock \emph{Mathematical Biosciences}, 28(3-4), 221--236.

\bibitem[{Liu et~al.(2019)Liu, Par\'{e}, Nedi\'{c}, Tang, Beck, and
  Ba\c{s}ar}]{bivirus_tac}
Liu, J., Par\'{e}, P.E., Nedi\'{c}, A., Tang, C.T., Beck, C.L., and Ba\c{s}ar,
  T. (2019).
\newblock Analysis and control of a continuous-time bi-virus model.
\newblock \emph{IEEE Transactions on Automatic Control}, 64(12), 4891--4906.

\bibitem[{Mai et~al.(2018)Mai, Battou, and Mills}]{suppress}
Mai, V.S., Battou, A., and Mills, K. (2018).
\newblock Distributed algorithm for suppressing epidemic spread in networks.
\newblock \emph{IEEE Control Systems Letters}, 2(3), 555--560.

\bibitem[{Mieghem et~al.(2009)Mieghem, Omic, and Kooij}]{van2009virus}
Mieghem, P.V., Omic, J., and Kooij, R. (2009).
\newblock Virus spread in networks.
\newblock \emph{IEEE/ACM Transactions on Networking}, 17(1), 1--14.

\bibitem[{Nowzari et~al.(2017)Nowzari, Preciado, and
  Pappas}]{preciado2017optimal}
Nowzari, C., Preciado, V.M., and Pappas, G.J. (2017).
\newblock Optimal resource allocation for control of networked epidemic models.
\newblock \emph{IEEE Transactions on Control of Network Systems}, 4(2),
  159--169.

\bibitem[{Par\'{e} et~al.(2015)Par\'{e}, Beck, and
  Nedi\'{c}}]{pare2015stability}
Par\'{e}, P.E., Beck, C.L., and Nedi\'{c}, A. (2015).
\newblock Stability analysis and control of virus spread over time--varying
  networks.
\newblock In \emph{Proc. 54th IEEE Conf. Decision Control}, 3554--3559.

\bibitem[{Par\'{e} et~al.(2018)Par\'{e}, Beck, and
  Nedi\'{c}}]{pare2016epidemic}
Par\'{e}, P.E., Beck, C.L., and Nedi\'{c}, A. (2018).
\newblock Epidemic processes over time-varying networks.
\newblock \emph{IEEE Transactions on Control of Network Systems}, 5(3),
  1322--1334.

\bibitem[{Par\'{e} et~al.(2020)Par\'{e}, Liu, Beck, Kirwan, and
  Ba\c{s}ar}]{pare2018dtjournal}
Par\'{e}, P.E., Liu, J., Beck, C.L., Kirwan, B.E., and Ba\c{s}ar, T. (2020).
\newblock Analysis, estimation, and validation of discrete-time epidemic
  processes.
\newblock \emph{IEEE Transactions on Control Systems Technology}, 28(1),
  79--93.

\bibitem[{Preciado et~al.(2014)Preciado, Zargham, Enyioha, Jadbabaie, and
  Pappas}]{PreciadoTCNS14}
Preciado, V.M., Zargham, M., Enyioha, C., Jadbabaie, A., and Pappas, G. (2014).
\newblock Optimal resource allocation for network protection against spreading
  processes.
\newblock \emph{IEEE Transaction on Control of Network Systems}, 1(1), 99--108.

\bibitem[{Santos et~al.(2015)Santos, Moura, and Xavier}]{santos2015bivirus}
Santos, A., Moura, J., and Xavier, J. (2015).
\newblock Bi-virus {SIS} epidemics over networks: qualitative analysis.
\newblock \emph{IEEE Transactions on Network Science and Engineering}, 2(1),
  17--29.

\bibitem[{Song et~al.(2018)Song, Hill, and Liu}]{song2018}
Song, Y., Hill, D.J., and Liu, T. (2018).
\newblock Network-based analysis of small-disturbance angle stability of power
  systems.
\newblock \emph{IEEE Transactions on Control of Network Systems}, 5(3),
  901--912.

\bibitem[{Torres and Roy(2018)}]{stu2}
Torres, J.A. and Roy, S. (2018).
\newblock Dominant eigenvalue minimization with trace preserving diagonal
  perturbation: Subset design problem.
\newblock \emph{Automatica}, 89, 160--168.

\bibitem[{Torres et~al.(2017)Torres, Roy, and Wan}]{stu1}
Torres, J.A., Roy, S., and Wan, Y. (2017).
\newblock Sparse resource allocation for linear network spread dynamics.
\newblock \emph{IEEE Transactions on Automatic Control}, 62(4), 1714--1728.

\bibitem[{Wan et~al.(2007)Wan, Roy, and Saberi}]{wan2007network}
Wan, Y., Roy, S., and Saberi, A. (2007).
\newblock Network design problems for controlling virus spread.
\newblock In \emph{Proc. 46th IEEE Conf. Decision Control}, 3925--3932.

\bibitem[{Wan et~al.(2008)Wan, Roy, and Saberi}]{wan2008designing}
Wan, Y., Roy, S., and Saberi, A. (2008).
\newblock Designing spatially heterogeneous strategies for control of virus
  spread.
\newblock \emph{IET Systems Biology}, 2(4), 184--201.

\bibitem[{Watkins et~al.(2018)Watkins, Nowzari, Preciado, and
  Pappas}]{watkins2016optimal}
Watkins, N.J., Nowzari, C., Preciado, V.M., and Pappas, G.J. (2018).
\newblock Optimal resource allocation for competitive spreading processes on
  bilayer networks.
\newblock \emph{IEEE Transactions on Control of Network Systems}, 5(1),
  298--307.

\bibitem[{Zelazo and B{\"u}rger(2014)}]{zelazo2014}
Zelazo, D. and B{\"u}rger, M. (2014).
\newblock On the definiteness of the weighted {L}aplacian and its connection to
  effective resistance.
\newblock In \emph{Proc. 53rd IEEE Conf. Decision Control}, 2895--2900.

\end{thebibliography}


\end{document}